\documentclass{article}
\usepackage {amsmath}
\usepackage {amsfonts}
\usepackage {amsthm}
\usepackage {indentfirst}
\usepackage {fullpage}

\title{An Exceptional Representation of $Sp(4,\Fq)$}
\author{Tanmay Deshpande}
\date{}

\newtheorem {thm} {Theorem} [section]
\newtheorem {prop} [thm] {Proposition}
\newtheorem {lem} [thm] {Lemma}
\newtheorem {cor} [thm] {Corollary}
\theoremstyle{definition}
\newtheorem {defn} [thm] {Definition}
\theoremstyle{remark}
\newtheorem {rk} [thm]  {Remark}

\newcommand{\f}{\mathbb}

\newcommand{\Fq} {\mathbb{F}_q}

\newcommand{\Fqcl} {\overline{\mathbb{F}}_q}

\newcommand{\F} {\mathbb{F}}

\newcommand{\C} {\mathbb{C}}

\newcommand{\Qlcl} {\overline{\mathbb{Q}}_l}
\newcommand{\beq}{\begin{equation}}
\newcommand{\eeq}{\end{equation}}
\newcommand{\bthm}{\begin {thm}}
\newcommand{\ethm}{\end {thm}}
\newcommand{\bprop}{\begin {prop}}
\newcommand{\eprop}{\end {prop}}
\newcommand{\bcor}{\begin {cor}}
\newcommand{\ecor}{\end {cor}}
\newcommand{\blem}{\begin{lem}}
\newcommand{\elem}{\end{lem}}
\newcommand{\bdefn}{\begin{defn}}
\newcommand{\edefn}{\end{defn}}
\newcommand{\brk}{\begin{rk}}
\newcommand{\erk}{\end{rk}}
\newcommand{\bpf}{\begin{proof}}
\newcommand{\epf}{\end{proof}}

\newcommand{\cU}{\cal{U}}
\newcommand{\cP}{\cal{P}}
\newcommand{\RTth}{R_{T,\theta}}
\newcommand{\RTthpr}{R_{T',\theta'}}

\newcommand{\into}{\hookrightarrow}

\begin{document}
\maketitle

\section{Introduction}
Let $p$ be an odd prime, and let $q$ be a power of $p$. Let $G = Sp(4)$. B.Srinivasan (in \cite{S68}) discovered an irreducible representation (denoted by $\theta_{10}$) of $Sp(4,\Fq)$ with the following remarkable combination of properties, namely it is cuspidal(Defn. \ref{cuspidal}), unipotent(Defn. \ref{unipotent}) as well as degenerate, i.e. it does not admit a Whittaker model(Defn. \ref{degenerate}). The groups $SL_n(\Fq)$ and $GL_n(\Fq)$ do not have any unipotent cuspidal representations and neither do they have any degenerate cuspidal representations. Hence the existence of such a representation for $Sp(4,\Fq)$ is somewhat surprising. 

We will describe a folklore construction of $\theta_{10}$, which is different from \cite{S68}. It is based on the Weil representation of $Sp(8, \Fq)$ and Howe duality. This article was a part of my master's thesis during my graduate studies at the University of Chicago. My advisor, V.Drinfeld, suggested that I publish this article in the e-print archive, since there were apparently no references for this construction of $\theta_{10}$.

\section{Construction}
Let $V$ be a four-dimensional symplectic vector space over $\Fq$ with a symplectic form $\langle \cdot, \cdot\rangle=\langle \cdot, \cdot\rangle_V$. Let $E={\f{F}}_{q^2}$ considered as a two-dimension vector space over $\Fq$ with a non-degenerate symmetric bilinear form corresponding to the norm, namely given by $\langle x,y \rangle_E = \frac{1}{2}(xy^q+x^qy)$. This is an anisotropic bilinear form. Moreover, any two-dimensional vector space over $\Fq$ with an anisotropic quadratic form is isomorphic to $E$. We see that the eight-dimensional space $V\otimes E$ inherits a natural symplectic form. Moreover, we have a natural map from $Sp(V) \times O(E)$ to $Sp(V\otimes E)$. Note that the group $SO(E)$ is just the cyclic group of order $q+1$, consisting of the norm 1 elements of $\f{F}_{q^2}^*$.

Let $\psi$ be a non-trivial character $\psi :\Fq \to \Qlcl^*$, where $l$ is a prime different from $p$. Let us now consider the corresponding Weil representation, $W$ of $Sp(V\otimes E)$. Thus we get an action of $Sp(V) \times O(E)$ on $W$. Let $L, L' \subset V$ be complementary Lagrangian subspaces. Then we can identify $W$ with the space of $\Qlcl$-valued functions on $L'\otimes E$.  Then for $t \in O(E)\into Sp(V\otimes E), v \in L'\otimes E$ and $f:L'\otimes E \to \Qlcl$, we have $(t \cdot f)(v) = f((1\otimes t^{-1})(v))$(See \ref{levi action}). For a character $\theta$ of $SO(E)$, let $W_\theta$ denote the $\theta$-isotypic part of $W$, i.e. $W_\theta = \{f|f((1\otimes t^{-1})(v)) = \theta(t)f(v) \mbox{ for all } t\in SO(E), v\in L'\otimes E\}$. Then $W_\theta$ is a representation of $Sp(V)$. We see that 
\beq\label{first} dim(W_\theta) = \left\{ \begin{array}{ll}
         \frac{q^4-1}{q+1}+1 = q(q^2-q+1) & \mbox{if $\theta = 1$}\\
         \frac{q^4-1}{q+1} = (q-1)(q^2+1) & \mbox{else}.\end{array} \right.\eeq 

Note that we have the element `conjugation', $(\sigma : \F_{q^2} \to \F_{q^2}) \in O(E)$. Let $\Lambda : W \to W$ denote the action of this element. Then $\Lambda$ commutes with the action of $Sp(V)$ and takes $W_\theta$ isomorphically to $W_{\theta^{-1}}$. Let $\nu$ be the quadratic character of $SO(E)$. Then for $\theta = 1$ or $\nu$, we have $\Lambda : W_\theta \to W_\theta$, and $\Lambda^2 = 1$. Let $W^{\pm}_\theta$ be the $(\pm1)$-eigenspace of $\Lambda|_{W_\theta}$. Let us now describe the dimensions of these four representations. 

\blem\label{dims}
$dim(W_\nu^{\pm})=\frac{1}{2}(q-1)(q^2+1)$, $dim(W_1^+)= \frac{1}{2}q(q^2+1)$ and $dim(W_1^-) = \frac{1}{2}q(q-1)^2$. 
\elem
\bpf
From the way $O(E)$ acts on $W$, we see that $(\Lambda \cdot f)(v) = f((1\otimes \sigma)(v))$, where $f\in W, v\in L'\otimes E$. Hence the trace of the operator $\Lambda = $ number of fixed points of $1\otimes \sigma: L'\otimes E \to L'\otimes E$. Now $\sigma$ is diagonalizable, with eigenvalues $\{1, -1\}$. Hence number of fixed points of $1\otimes \sigma$ $=$ number of elements in the $1$-eigenspace of $1\otimes \sigma = q^2$. Hence the trace of $\Lambda = q^2$. Similarly, $W_1$ can be identified with the space of $\Qlcl$-valued functions on $(L'\otimes E)/SO(E)$, and as before the trace of $\Lambda|_{W_1} = $ number of fixed points of $1\otimes \sigma: (L'\otimes E)/SO(E) \to (L'\otimes E)/SO(E)=|S/SO(E)|$, where the set $S=\{v\in L'\otimes E|(1\otimes \sigma)(v)=(1\otimes \lambda)(v) \hbox{ for some } \lambda \in SO(E)\}$ i.e. $S/SO(E)$ is precisely the set of self-conjugate orbits. Let $(v'_1, v'_2)$ be a basis of $L'$. Then any vector $v \in L'\otimes E$ can be uniquely written as $v=v'_1\otimes e_1 + v'_2\otimes e_2$, where $e_1, e_2 \in E$ are uniquely determined by $v$. Then such a $v \in S \iff v'_1\otimes e_1^q + v'_2\otimes e_2^q = v'_1\otimes \lambda e_1 + v'_2\otimes \lambda e_2$ for some $\lambda \in SO(E) \iff e_1^{q-1} = e_2^{q-1}$ or one of $e_1, e_2 = 0$ $ \iff e_1, e_2$ linearly dependent $\iff v$ is decomposable, i.e. $v = v'\otimes e$ for some $v' \in L', e \in E$. Hence $S$ is precisely the set of decomposable vectors. Now the number of non-zero decomposable vectors in $L'\otimes E = \frac{(q^2-1)^2}{q-1}$ and hence $|(S\backslash\{0\})/SO(E)|= q^2-1$. So we see that $|S/SO(E)|= q^2$ , and hence trace of $\Lambda|_{W_1} = q^2$ as well. Hence from the way $\Lambda$ acts on $W$, we see that trace of  $\Lambda|_{W_\nu}=0$. The lemma now follows from (\ref{first}), since the trace of $\Lambda|_{W_\theta} = dim(W_\theta^+) - dim(W_\theta^-)$ for $\theta = 1$ or $\nu$.
\epf

\bprop\label{irreducible}
The representations $W_\theta\cong W_{\theta^{-1}}$ for $\theta \neq 1, \nu$, $W_1^{\pm}$ and $W_\nu^{\pm}$ are all irreducible and distinct.
\eprop
\bpf
We see that as an $Sp(V)$ representation we have
\beq
W = 2.\left(\bigoplus_{\langle \theta\rangle, \theta \neq 1, \nu}{W_\theta}\right) \oplus W^+_1 \oplus W^-_1 \oplus W^+_\nu \oplus W^-_\nu.
\eeq
Hence we must have $\langle W,W\rangle_{Sp(V)} \geq 4(q-1)/2+4=2q+2$, since all the summands above are non-zero. On the other hand, we now show that $dim(End_{Sp(V)}(W)) = 2q+2$. Let $\mathcal{A}$ be the group algebra of the Heisenberg group. Let $\mathcal{A}_\psi$ denote the quotient of $\mathcal{A}$ corresponding to the central character $\psi$. Since $W$ is the space of the irreducible representation of the Heisenberg group with central character $\psi$, we get a canonical isomorphism $\mathcal{A}_\psi \tilde{\rightarrow} End_{\C}(W)$ that is $Sp(V\otimes E)$-equivariant. Now as a representation of $Sp(V\otimes E)$, $\mathcal{A}_\psi$ identifies with the space of $\Qlcl$-valued functions on $V\otimes E$. Now $End_{Sp(V)}(W) = End_{\C}(W)^{Sp(V)}$. Hence we see that $dim(End_{Sp(V)}(W)) = dim((End_{\C}(W)^{Sp(V)}) =dim(\mathcal{A}_\psi^{Sp(V)}) = $ number of $Sp(V)$-orbits in $V\otimes E$. Let $(e_1,e_2)$ be a basis of $E$. Then as before, an elements of $V\otimes E$ can be uniquely written in the form $v_1\otimes e_1+v_2\otimes e_2$ with $v_1, v_2 \in V$. Then, we have the orbit $\{0\}$. The set of orbits of non-zero decomposable vectors can be identified with $\f{P}(E)$, and finally the set of orbits of indecomposable elements can be identified with $\Fq$ via the correspondence $Sp(V)\cdot(v_1\otimes e_1+v_2\otimes e_2){\leftrightarrow}\langle v_1, v_2\rangle_V$. So we see that the number of orbits is exactly $2q+2$ i.e. $\langle W,W \rangle_{Sp(V)} = 2q+2$. Hence we conclude that all the summands in the decomposition above must be irreducible and distinct. 
\epf
In particular, $W_1^-$ is an irreducible $d=$\mbox{$\frac{1}{2}q(q-1)^2$-dimensional} representation of $Sp(V)$. Prop.\ref{little groups} gives another proof of the irreducibility of this representation. For historical reasons, let us denote its character by $\theta_{10}$. We will prove that this representation is cuspidal, degenerate and unipotent. Let us first study the space $W^-_1$. Let $S'$ be the set of indecomposable vectors in $L'\otimes E$. As we have seen in the proof of \ref{dims}, $S'/SO(E) \subset (L'\otimes E)/SO(E)$ is precisely the set of $SO(E)$-orbits that are not self-conjugate. Let $\mathcal{O}_1, \mathcal{O}'_1, \mathcal{O}_2, \mathcal{O}'_2, \cdots, \mathcal{O}_d, \mathcal{O}'_d$ be all such orbits, where $\mathcal{O}'_i = (1\otimes \sigma)\mathcal{O}_i$. Then it is clear that the functions $\delta_i=\delta_{\mathcal{O}_i} - \delta_{\mathcal{O}'_i}$ form a basis of $W^-_1$, where for $X\subset L'\otimes E$, $\delta_X$ denotes the function that takes the value $1$ on $X$ and $0$ elsewhere.

\section{Parabolic Subgroups of $Sp(V)$}
The Weyl group of $Sp(4)$ is isomorphic to the dihedral group $D_8$.
Let $0 \subset L_1 \subset L(=L^{\perp}) \subset L_1^{\perp} \subset V$ be a complete flag in $V$ and let $L'$ be a complementary Lagrangian subspace to the Lagrangian subspace $L \subset V$. Now the stabilizer $B_0$ of this complete flag is a Borel subgroup. Let $U_0$ be its unipotent radical. Then $U_0$ is a maximal unipotent subgroup and its order is $q^4$. Let $(v_1,v_2)$ be a basis of $L$ such that $v_1\in L_1$, and let $(v_3,v_4)$ be a basis of $L'$ such that the matrix of $\langle,\rangle_V$ with respect to the basis $(v_1,v_2,v_3,v_4)$ is \( \left( \begin{smallmatrix}
0 & 0 & 0 & 1\\
0 & 0 & 1 & 0\\
0 &-1 & 0 & 0\\
-1& 0 & 0 & 0 \end{smallmatrix} \right).\) Thus we have identified $Sp(V)$ with $Sp(4, \Fq)$. We see $v_1 \in L_1, v_2 \in L, v_3 \in L_1^\perp, v_4 \in V$. Then with respect to this basis, $U_0$ is the group of matrices of the type \( \left( \begin{smallmatrix}
1 & -\alpha & \beta & \mu\\
0 & 1 & \lambda & \lambda\alpha+\beta\\
0 & 0 & 1 & \alpha\\
0& 0 & 0 & 1 \end{smallmatrix} \right)\) where $(\lambda, \alpha, \mu, \beta) \in \Fq^4$. Let $T_0$ be the torus of diagonal matrices \( \left( \begin{smallmatrix}
a & 0 & 0 & 0\\
0 & b & 0 & 0\\
0 & 0 & b^{-1} & 0\\
0& 0 & 0 & a^{-1} \end{smallmatrix} \right)\) in $Sp(V)$ with respect to this basis. Then we can identify\footnote{Here we use the convention that for certain nice $X\subset Sp(V)$, $\underline{X}$ denotes the obvious $F$-stable subvariety of $Sp(4)$ such that $\underline{X}^F=X$.} $\underline{T_0}$ with $\f{G}_m \times \f{G}_m$, and hence we can identify $Hom(\underline{T_0}, \f{G}_m)$ with $\f{Z}^2$. With this identification, the roots of $Sp(4)$ are $\{\pm(1,-1),\pm(0,2),\pm(1,1),\pm(2,0)\}$, and the choice of positive roots (implicit in this notation) is forced by our choice of the Borel subgroup. The simple roots are $r_1=(1,-1)$ and $r_2=(0,2)$. The other positive roots are $r_3 = r_1+r_2$ and $r_4 = 2r_1+r_2$. The Weyl group $W(\underline{T_0})$ is generated by the two simple reflections $s_1 = $\( \left( \begin{smallmatrix}
0 & 1 & 0 & 0\\
1 & 0 & 0 & 0\\
0 & 0 & 0 & 1\\
0 & 0 & 1 & 0 \end{smallmatrix} \right)\)and $s_2 = $ \( \left( \begin{smallmatrix}
1 & 0 & 0 & 0\\
0 & 0 & 1 & 0\\
0 &-1 & 0 & 0\\
0 & 0 & 0 & 1 \end{smallmatrix} \right)\) corresponding to the simple roots $r_1, r_2$ respectively. Let us now describe the various root subgroups. $$\underline{U_{r_1}}=\left\{\left( \begin{smallmatrix}
1 & -\alpha & 0 & 0\\
0 & 1 & 0 & 0\\
0 & 0 & 1 & \alpha\\
0& 0 & 0 & 1 \end{smallmatrix} \right)\right\}, \underline{U_{-r_1}} = \underline{U_{r_1}^{T}}.$$
$$\underline{U_{r_2}}=\left\{\left( \begin{smallmatrix}
1 & 0 & 0 & 0\\
0 & 1 & \lambda & 0\\
0 & 0 & 1 & 0\\
0 & 0 & 0 & 1 \end{smallmatrix} \right)\right\}, \underline{U_{-r_2}} = \underline{U_{r_2}^{T}}.$$
$$\underline{U_{r_3}}=\left\{\left( \begin{smallmatrix}
1 & 0 & \beta & 0\\
0 & 1 & 0 & \beta\\
0 & 0 & 1 & 0\\
0 & 0 & 0 & 1 \end{smallmatrix} \right)\right\}, \underline{U_{-r_3}} = \underline{U_{r_3}^{T}}.$$
$$\underline{U_{r_4}}=\left\{\left( \begin{smallmatrix}
1 & 0 & 0 & \mu\\
0 & 1 & 0 & 0\\
0 & 0 & 1 & 0\\
0 & 0 & 0 & 1 \end{smallmatrix} \right)\right\}, \underline{U_{-r_4}} = \underline{U_{r_4}^{T}}.$$

Let $P_1$ be the Siegel parabolic subgroup corresponding to the Lagrangian subspace $L$, i.e. the subgroup of elements of $Sp(V)$ that leave the sub-flag  $0 \subset L \subset V$ invariant. This is the parabolic subgroup $B_0\cup B_0s_1B_0$. The unipotent radical $U_1$ of $P_1$ consists of those elements of $Sp(V)$ that act as identity on $L$. We have $U_1=U_{r_2}U_{r_3}U_{r_4}$. It consists of the matrices \( \left( \begin{smallmatrix}
1 & 0 & \beta & \mu\\
0 & 1 & \lambda &\beta\\
0 & 0 & 1 & 0\\
0& 0 & 0 & 1 \end{smallmatrix} \right)\).
\blem
$U_1$ can be naturally identified with the additive group of symmetric bilinear forms on $L'$. 
\elem
\bpf
For $g \in Sp(V)$, we have $ker(g-1)^{\perp} = Im(g-1)$. Suppose $g \in U_1$. Then $L \subset Ker(g-1)$ and $Im(g-1) \subset L$. Thus $g-1$ induces a map $A_g:L'\cong V/L \to L$. On the other hand, a bilinear form $A$ on $L'=V/L$ given by a map $A:V/L \to L$, induces a map $A':V \to V$. Then we have that $\langle (1+A')v, (1+A')w\rangle = \langle v, w\rangle \Longleftrightarrow \langle v, A'w\rangle + \langle A'v, w\rangle = 0$. Hence we see that $1+A' \in Sp(V) \Longleftrightarrow A$ is a symmetric bilinear form on $L'$. Moreover, if $A_1, A_2$ are two symmetric bilinear forms on $L'$, then $(1+A'_1)(1+A'_2) = 1+(A_1+A_2)'$. Hence we have identified the group $U_1$, with the group of symmetric bilinear forms on $L'$.
\epf
So we see that in fact, $U_1$ is a 3-dimensional vector space over $\Fq$. For $g \in U_1$, let $\langle \cdot ,\cdot \rangle_g$ denote the corresponding bilinear form on $L'$.

Let $P_2$ be the stabilizer of the flag $0 \subset L_1 \subset L_1^{\perp} \subset V$. This is the parabolic subgroup $B_0\cup B_0s_2B_0$. Let $U_2$ be its unipotent radical. Then $U_2=${$\{g \in Sp(V)|(g-1)L_1=0 \hbox{ and } (g-1)L_1^{\perp} \subset L_1\}$.} We have $|U_2|=q^3$. We have $U_2=U_{r_1}U_{r_3}U_{r_4}$. We see that $U_2$ consists of the matrices \( \left( \begin{smallmatrix}
1 & -\alpha & \beta & \mu\\
0 & 1 & 0 & \beta\\
0 & 0 & 1 & \alpha\\
0 & 0 & 0 & 1 \end{smallmatrix} \right)\).

Now $B_0, P_1$ and $P_2$ are all the proper parabolic subgroups containing $B_0$.

Let $U'_0 = U_1 \cap U_2$. Then $U'_0$ is the commutator subgroup of $U_0$. For $g \in U_1$, we see that $g \in U_2 \Longleftrightarrow \langle L_1^{\perp}, (g-1)L_1^{\perp}\rangle = 0 \Longleftrightarrow \langle L_1^{\perp}/L, L_1^{\perp}/L \rangle_g = 0$. Hence $U'_0$ can be identified with the group of all symmetric bilinear forms on $L'$ such that $L_1^{\perp}/L \subset L'$ is an isotropic subspace with respect to that form. We have $|U'_0| = q^2$. We have $U'_0=U_{r_3}U_{r_4}$. It consists of the matrices \( \left( \begin{smallmatrix}
1 & 0 & \beta & \mu\\
0 & 1 & 0 & \beta\\
0 & 0 & 1 & 0\\
0& 0 & 0 & 1 \end{smallmatrix} \right)\). Let $U''_0 \subset U'_0$ be the group of all symmetric bilinear forms that contain $L_1^{\perp}/L$ in their kernels. Then in fact $U''_0$ is the center of $U_0$. We have $U''_0 = U_{r_4}\cong \Fq$. It consists of the matrices \( \left( \begin{smallmatrix}
1 & 0 & 0 & \mu\\
0 & 1 & 0 & 0\\
0 & 0 & 1 & 0\\
0& 0 & 0 & 1 \end{smallmatrix} \right)\).

\section{Cuspidality and Degeneracy}
Let us first recall the definitions of cuspidality and degeneracy for finite groups of Lie type.
\bdefn\label{cuspidal}
Let $G$ be a connected reductive group over $\Fqcl$ that has an $\Fq$-structure given by a geometric Frobenius endomorphism $F$. Let $\rho$ be a representation of $G^F$. We say that $\rho$ is cuspidal, if for any proper $F$-stable parabolic subgroup $P$ of $G$ with unipotent radical $U$, $\rho|_{U^F}$ does not contain the trivial representation of $U^F$.
\edefn

Let us now recall the definition of a Whittaker model. For simplicity, let us assume that the group $G$ is split over $\Fq$. Let $U_0$ be a maximal $F$-stable unipotent subgroup of $G$ contained in an $F$-stable Borel subgroup $B_0$. Let $B_0 = T_0U_0$, where $T_0$ is $F$-stable. The group $U_0/[U_0, U_0]$ is isomorphic to a sum of copies of the groups $\f{G}_a$, one copy for each simple root of $G$ with respect to the pair $(T_0, B_0)$. $G$ being split implies that these copies of $\f{G}_a$ are $F$-stable. Let $\xi: U_0^F/[U_0, U_0]^F \to \Qlcl^*$ be a non-degenerate character, i.e. the restriction of $\xi$ to each copy of $\f{G}_a^F$ is non-trivial. We can consider $\xi$ as a character of $U_0^F$. We call such characters $\xi$ of $U_0^F$ as non-degenerate characters. Then the representation $Ind_{U_0^F}^{G^F}(\xi)$ is multiplicity free.(See \cite{St}, p. 258-262.)
\bdefn\label{degenerate}
Let $\rho$ be an irreducible representation of $G^F$. We say that $\rho$ admits a Whittaker model, if $\rho|_{U_0^F}$ contains a one-dimensional non-degenerate character $\xi$ of $U_0^F$. We say that $\rho$ is degenerate if it does not admit a Whittaker model.
\edefn
\brk
By Frobenius reciprocity and the fact that $Ind_{U_0^F}^{G^F}(\xi)$ is multiplicity free, we see that if $\rho_{U_0}$ contains such a non-degenerate $\xi$, then we must have $\langle \xi, \rho|_{U^F_0}\rangle_{U_0^F} = 1$. 
\erk
\brk\label{remark}
If $\rho$ is such that $\rho|_{U_0^F}$ does not contain any one-dimensional characters, then $\rho$ must be degenerate.
\erk

Let us now return to the group $Sp(V)$. We now study how $U_1$ acts on the Weil representation $W$. Now $L\otimes E$ and $L'\otimes E$ are complementary Lagrangian subspaces of the symplectic vector space $V\otimes E$. So, as before, we can identify the space $W$ with \mbox{$\{f|f : L'\otimes E \to \Qlcl \}$.} Let $\cal{P}$ be the Seigel parabolic subgroup of $Sp(V\otimes E)$ corresponding to $L\otimes E$ i.e. the stabilizer of this Lagrangian subspace. Let $\cal{U}$ be its unipotent radical. Then exactly as before, we may identify $\cal{U}$ with the additive group of symmetric bilinear forms on $L'\otimes E$. We have the inclusion $Sp(V)\hookrightarrow Sp(V\otimes E)$  given by $g \mapsto g\otimes 1$. This induces the inclusions $P_1 \hookrightarrow \cP$ and $U_1 \hookrightarrow \cU$. In terms of bilinear forms, we have $\langle,\rangle_{g\otimes 1} = \langle,\rangle_g\otimes \langle,\rangle_E$ for $g \in U_1$. We now recall how $\cU$ acts on $W$.
\bprop[See \cite{G}, 2.8.]\label{action}
Let $g \in \cU$ and $f:L'\otimes E \to \Qlcl$. Then
\beq
(g\cdot f)(v) = \psi\left(\frac{\langle v,v\rangle_g}{2}\right)f(v).
\eeq
\eprop

The cuspidality and degeneracy of $W_1^-$ follow immediately from the following:
\blem
The restriction $W_1^-|_{U'_0}$ does not contain the trivial representation of $U'_0$, i.e. $W_1^-$ has no non-zero $U'_0$-fixed vector. In fact, $W^-_1$ does not even have any non-zero $U''_0$-fixed vectors.
\elem 
\begin{proof}
Let $f \in W$ be a fixed vector for the action $U''_0$. Suppose $v \in L'\otimes E$ is such that $f(v) \neq 0$. Then from \ref{action} we conclude that $\psi\left(\frac{\langle v,v\rangle_{g\otimes 1}}{2}\right) = 1$ for all $g \in U''_0$. Hence for all $\alpha \in \Fq$, $$\psi\left(\frac{\langle v,v\rangle_{\alpha\cdot g\otimes 1}}{2}\right) = \psi\left(\alpha\cdot \left(\frac{\langle v,v\rangle_{g\otimes 1}}{2}\right)\right) = 1.$$ Hence we must have $\langle v,v\rangle_{g\otimes 1} = 0$ for all $g \in U''_0$. So we see that $v$ must in fact lie in $(L_1^{\perp}/L)\otimes E$, say $v = v'\otimes e$, where $e \in E=\F_{q^2}$. Then $e^{q-1}$ is a norm 1 element of $\F_{q^2}$, i.e. an element of $SO(E)$. Now suppose $f \in W^-_1$. Then $f((1\otimes e^{q-1})(v'\otimes e))=f(v'\otimes e)$ since $f\in W_1$. Hence $f(v'\otimes e^q) = f(v)$. On the other hand, $\Lambda f = -f$, hence $f((1\otimes \sigma)(v'\otimes e))=-f(v'\otimes e)$, i.e. $f(v'\otimes e^q) = -f(v)$. Hence we arrive at a contradiction.
\end{proof}

\bprop
The representation $W_1^-$ is cuspidal and degenerate.
\eprop
\bpf
The lemma implies that none of the unipotent radicals $U_0,U_1,U_2$ has a fixed vector in $W^-_1$ , or in other
words $W^-_1$ is cuspidal. On the other hand, $W^-_1|_{U_0}$  cannot have a one-dimensional $U_0$-invariant subspace, for
otherwise, the commutator $U'_0$ would have to fix some vector. Hence by \ref{remark}, we conclude that the representation is
degenerate.
\epf

Next, we show that in fact the representation of $P_1$ on $W_1^-$ is irreducible. For this, let us study the restriction of $W_1^-$ to $U_1$. We will need the following:
\blem\label{easy}
Let $v, w \in S'$ i.e. $v,w$ indecomposable. Then $\langle v, v \rangle_{g\otimes 1} = \langle w, w \rangle_{g\otimes 1}$ for all $g \in U_1 \iff v,w$ are in the same $O(E)$-orbit. For an $O(E)$-orbit $\mathcal{O}$ in $S'$, let $\phi_{\mathcal{O}}:U_1 \to \Qlcl^*$ be the character $g\mapsto \psi\left(\frac{\langle v, v \rangle_{g\otimes 1}}{2}\right)$ where $v\in \mathcal{O}$.
\elem
\bpf
Let $v = v_3\otimes e_1 + v_4\otimes e_2$, $w = v_3\otimes f_1 + v_4\otimes f_2$ where $v_3, v_4$ are as before. Then since we have $\langle,\rangle_{g\otimes 1} = \langle,\rangle_g\otimes \langle,\rangle_E$ for all $g \in U_1$, we conclude from the hypothesis of the lemma, that we must have $\langle e_i,e_j\rangle_E = \langle f_i,f_j\rangle_E$. Now $v$ indecomposable$\implies e_1,e_2$ linearly independent. Similarly $f_1,f_2$ also linearly independent. Hence we conclude that there exists $t \in O(E)$ such that $(1\otimes t)(v) = w$. The converse is obvious. 
\epf

\blem\label{stupid}
Let $\phi:U_1 \to \Qlcl^*$ be a character of $U_1$. Let $f \in W$. Then $U_1$ acts on $f$ by $\phi \iff $ for all $v \in L'\otimes E$ such that $f(v)\neq 0$ we have $\psi\left(\frac{\langle v,v \rangle_{g\otimes 1}}{2}\right) = \phi(g)$. Hence $\Qlcl\cdot f=\langle f\rangle \subset W$ is $U_1$-invariant $\iff$ for all $v,w \in L'\otimes E$ where $f$ does not vanish, we have  $\langle v,v \rangle_{g\otimes 1}=\langle w,w \rangle_{g\otimes 1}$ for all $g \in U_1$.
\elem
\bpf
This is immediate from \ref{action}.
\epf

\bprop
We have the decomposition $W_1^-|_{U_1} = \bigoplus_{i}{\langle \delta_i \rangle}$ as $U_1$-modules. Let $\hat{\mathcal{O}}_i \subset S'$ be the $O(E)$-orbit ${\mathcal{O}}_i\cup {\mathcal{O}}'_i$. Then $\langle \delta_i \rangle \cong \phi_{\hat{\mathcal{O}}_i}$ are distinct as $U_1$-modules.
\eprop
\bpf
That $\langle \delta_i \rangle$ is $U_1$-invariant follows from \ref{stupid} and one direction of \ref{easy}. The second assertion in the proposition follows  from the other direction of \ref{easy}.
\epf

Let $M_1$ be the stabilizer in $Sp(V)$ of $(L, L')$, and let $\mathcal{M} \subset Sp(V\otimes E)$ be the stabilizer of $(L\otimes E, L'\otimes E)$. Then $M_1$ is a Levi subgroup of $P_1$, while $\mathcal{M}$ is a Levi subgroup of $\mathcal{P}$. So we have $P_1 = M_1U_1$ and $\mathcal{P}=\mathcal{M}\mathcal{U}$. $M_1$ can be identified with $GL(L)$ or $GL(L')$, and similarly $\mathcal{M}$ can be identified with $GL(L\otimes E)$ or $GL(L'\otimes E)$. We have $M_1 \into \mathcal{M}$. We now state how $\mathcal{M}$ acts on $W$.
\bprop[See \cite{G}, 2.7.]\label{levi action}
Let $g \in \mathcal{M}$ be thought of as an element of $GL(L'\otimes E)$. Let $f \in W$. Then for $v \in L'\otimes E$ we have
\beq
(g\cdot f)(v) = \chi(det(g))f(g^{-1}v),
\eeq
where $\chi$ is the quadratic character of $\Fq^*$.\\If $g \in M_1$, we have $det(g\otimes 1) = det(g)^2$, hence
\beq
(g\cdot f)(v) = f((g^{-1}\otimes 1)v).
\eeq

\eprop
It is easy to see that $M_1$ acts simply transitively on $S' \subset L'\otimes E$. So $M_1$ acts transitively on $S'/O(E)$, which we can identify with $\{1,2,\cdots,d\}$, and the stabilizer of the element $1$, i.e. of $\hat{\mathcal{O}}_1\in S'/O(E)$ in $M_1$ is a subgroup $O$ isomorphic to $O(E)$. From above, we see that $g \in M_1$ takes the space $\langle\delta_i \rangle$ to the space $\langle\delta_{g\cdot i}\rangle$. We now apply the little groups method(See \cite{S}, section 8.2) in our setting to get the following result.
\bprop\label{little groups}
The restriction $W_1^-|_{P_1}$ is an irreducible representation of $P_1$. It is induced from a one-dimensional character of a certain subgroup of $P_1$.
\eprop
\brk
This gives another proof of the irreducibility of $W_1^-$.
\erk
\bpf
From what we have proven so far, we see that the characters $\phi_{\hat{\mathcal{O}}_i}(\cong\langle \delta_i\rangle$ as $U_1$-modules) of $U_1$ form a single $M_1$-orbit. The stabilizer of $\phi=\phi_{\hat{\mathcal{O}}_1}$ in $M_1$ is the subgroup $O$ isomorphic to $O(E)$. Then from \ref{levi action}, we see that $O$ acts on $\langle \delta_1 \rangle$ by the sign representation $\epsilon: O\to \{\pm 1\}$. Hence the subgroup $OU_1$ is the stabilizer in $P_1$ of the one-dimensional subspace $\langle \delta_1 \rangle$. The character $\phi'$, or the action of $OU_1$ on this subspace is given by $\phi'(hg)=\epsilon(h)\phi(g)$ for $h \in O, g\in U_1$. Then we have $W_1^-|_{P_1} \cong Ind^{P_1}_{OU_1}(\phi')$ and that it is irreducible.
\epf

\section{Deligne-Lusztig Theory}
Let us recall some results from Deligne-Lusztig theory that are relavent. Let $G$ be a connected reductive group over $\Fqcl$ provided with an $\Fq$-structure given by a geometric Frobenius morphism $F:G\to G$. Let $(T_0, B_0)$ be a pair consisting of an $F$-stable maximal torus and an $F$-stable Borel subgroup containing it respectively. Let $B_0 = T_0U_0$. For an $F$-stable maximal torus $T$, let $\hat{T}^F$ denote the group of characters of $T^F$ with values in $\Qlcl$. Let $T, T'$ be two $F$-stable maximal tori. We define $N(T,T') = \{g \in G| g^{-1}Tg=T'\}$. Define $W(T,T')=T\backslash N(T, T') = N(T,T')/T'$. Note that since $T$, $T'$ are $F$-stable, $N(T,T')$ will also be $F$-stable and we will have an induced action of $F$ on $W(T,T')$. Then using Lang's Theorem, we observe that $W(T,T')^F$ can be identified with $T^F\backslash N(T,T')^F$ or with $N(T,T')^F/T'^F$.

For each integer $n>0$, we have the norm map $N_n: T^{F^n} \to T^F$. 
\bdefn[\cite{DL}, Defn. 5.5.]
Let $\theta, \theta'$ be characters of $T^F, T'^F$ respectively. We say that $(T,\theta)$ and $(T', \theta')$ are geometrically conjugate if there exists an integer $n>0$ and $g \in G^{F^n}$ such that ${}^{g}T' = T$ and ${}^{g}(\theta' \circ N_n) = \theta \circ N_n$.
\edefn
\bthm[\cite{DL}, Cor. 6.3.]\label{strong orthogonality}
Let $\theta, \theta'$ be characters of $T^F, T'^F$ respectively. If $(T,\theta), (T', \theta')$ are not geometrically conjugate, then the virtual representations $R_{T,\theta}$ and $R_{T',\theta'}$ are disjoint, i.e. have no irreducible components in common.
\ethm

\bthm[\cite{DL}, Thm. 6.8.] \label{weak orthogonality}
\beq
	\langle R_{T,\theta}, R_{T',\theta'}\rangle=|\{w\in W(T,T')^F|{}^{w}\theta'=\theta \}|.
\eeq
In particular, $\langle \RTth,\RTthpr\rangle=0$ if $(T, \theta)$, $(T', \theta')$ are not $G^F$-conjugate.
\ethm

\brk
This does not mean that $R_{T,\theta}$ and $\RTthpr$ are disjoint, since $\RTth, \RTthpr$ are only virtual characters.
\erk

\bdefn
	We say that a character $\theta$ of ${T}^F$ is in general position, or is regular if it is not fixed by any non-trivial element of $W(T,T)^F$.
\edefn

\bcor
If $\theta$ is regular then $\pm \RTth$ is irreducible.
\ecor

Let $St_G$ denote the Steinberg representation of $G^F$. For an $F$-stable torus $T$, let $\epsilon_T = (-1)^s$, where $s = $ the dimension of the split part of $T$. We let $\epsilon_G = \epsilon_{T_0}$.

\bthm[\cite{DL}, Thm. 7.1.]\label{dimension}
\beq
dim(\RTth)=Q_T(1)=\epsilon_G\epsilon_T\frac{|G^F|}{|U^F_0||T^F|}=\epsilon_G\epsilon_T\frac{|G^F|}{St_G(1)|T^F|}
\eeq
\ethm
Let us recall the definition of unipotence.
\bdefn\label{unipotent}
Let $\rho$ be an irreducible representation of the group $G^F$. We say that $\rho$ is \emph{unipotent} if it occurs in some virtual character $R_{T, 1}$ for some $F$-stable maximal torus $T$.
\edefn
\brk
By \ref{strong orthogonality} we see that a unipotent representation cannot occur in $R_{T, \theta}$ if $\theta \neq 1$.
\erk

We will make use of the following result to show that $\theta_{10}$ is unipotent.
\bprop[\cite{DL}, Cor. 7.6.]\label{fairy tale}
Let $\rho$ be any virtual character of $G^F$ and let $s \in G^F$ be semisimple. Then
\beq
\rho(s)=\frac{1}{St_G(s)}\sum_{T\ni s}{\sum_{\theta \in \hat{T}^F}{\epsilon_G\epsilon_T}\theta(s)\langle \rho, R_{T,\theta}\rangle}.
\eeq
In particular, if $s$ is regular semisimple and if $T$ is the unique maximal torus containing it, then
\beq
\rho(s)= \sum_{\theta \in \hat{T}^F}{\theta(s)\langle \rho, R_{T,\theta} \rangle}.
\eeq
\eprop

\section{Unipotence}
Let us now return to the case where $G=Sp(4)$. 
We will make use of the following formula for the values of the character $\eta$ of the Weil representation on a certain subset of $Sp(V\otimes E)$. This formula was obtained by S.Gurevich and R.Hadani as a consequence of their algebro-geometric approach to the Weil representation(See \cite{GH}).

\bprop[See \cite{GH}; \cite{T}, Rem. 1.3.]\label{character formula}
Let $g \in Sp(V\otimes E)$ be such that $g-1$ is invertible. Let $\chi$ be the quadratic character of $\Fq^*$. Then
\beq
\eta(g) = \chi(det(g-1)).
\eeq 
\eprop
We will make use of this formula to compute $\theta_{10}(s)$, where $s$ is any regular element of a certain $F$-stable maximal torus $T$. Let $\kappa$ be a generator of $\F_{q^4}^*$ and let $\zeta=\kappa^{q^2-1}$. Then let $T$ be an $F$-stable maximal torus in $Sp(4)$ such that $T^F \cong \langle \zeta \rangle$.(See \cite{S68}, 3.2.) The group $T^F$ is conjugate in $Sp(4,\Fqcl)$ to the subgroup $H\subset Sp(4, \Fqcl)$ generated by the matrix \( \left( \begin{smallmatrix}
\zeta^q & 0 & 0 & 0\\
0 & \zeta & 0 & 0\\
0 & 0 & \zeta^{-1} & 0\\
0 & 0 & 0 & \zeta^{-q} \end{smallmatrix} \right) \in Sp(4,\Fqcl)\). $N(T,T)^F$ is conjugate to the subgroup of $Sp(4, \Fqcl)$ generated by $H$ and \( \left( \begin{smallmatrix}
0 & 1 & 0 & 0\\
0 & 0 & 0 & 1\\
-1& 0 & 0 & 0\\
0 & 0 & 1 & 0 \end{smallmatrix} \right)\). $W(T,T)^F$ is the cyclic group of order $4$, and its generator acts on $T^F$ and $\hat{T}^F$ by taking $q$-th powers. Let $s \in T^F$ correspond to $\gamma \in \langle\zeta\rangle$. Then by the description of $T^F$ above, it follows that the eigenvalues of $s:V \to V$ are $\{\gamma, \gamma^q, \gamma^{q^2}=\gamma^{-1}, \gamma^{q^3}=\gamma^{-q}\}$. 

\bprop\label{rss}
Let $s \in T^F, s \neq \pm1$. Then $s$ is regular semisimple and $\theta_{10}(s)=1$.
\eprop
\bpf
$W_1^-$ was defined as the $\epsilon$-isotypic component of $W$, where $\epsilon:O(E) \to \{\pm1\}$ is the `sign' character. Hence we see that for $g \in Sp(V)$ we have
\beq\label{eqn}
\theta_{10}(g) = \frac{1}{2q+2}\sum_{t\in O(E)}{\epsilon(t^{-1})\eta(g\otimes t)} = \frac{1}{2q+2}\left(\sum_{t\in SO(E)}{\eta(g\otimes t)}-\sum_{t\in O(E)\backslash SO(E)}{\eta(g\otimes t)}\right).
\eeq
In view of this formula, \ref{rss} would immediately follow once we have the following:
\blem
Let $s$ be as in the Proposition \ref{rss}. Then we have
\beq
\eta(s\otimes t) = \left\{ \begin{array}{ll}
         1 & \mbox{if $t \in SO(E)$}\\
        -1 & \mbox{if $t \in O(E)\backslash SO(E)$}.\end{array} \right.
\eeq
\elem

\bpf
For $s$ as above, and $t \in SO(E)$ (having eigenvalues $\{t, t^q=t^{-1}\}$), $s\otimes t - 1$ is invertible (since $\gamma\neq \pm1$) with eigenvalues, \mbox{$\{t\gamma-1, t\gamma^q-1, t\gamma^{-1}-1, t\gamma^{-q}-1, t^{-1}\gamma-1, t^{-1}\gamma^q-1,  t^{-1}\gamma^{-1}-1, t^{-1}\gamma^{-q}-1\}$}. $$det(s\otimes t-1) = \frac{(t\gamma-1)^2}{t\gamma}\frac{(t\gamma^q-1)^2}{t\gamma^q}\frac{(t\gamma^{-1}-1)^2}{t\gamma^{-1}}\frac{(t\gamma^{-q}-1)^2}{t\gamma^{-q}} = \left(\frac{(t\gamma-1)(t\gamma^q-1)(t\gamma^{-1}-1)(t\gamma^{-q}-1)}{t^2}\right)^2.$$ Now $$\frac{(t\gamma-1)(t\gamma^q-1)(t\gamma^{-1}-1)(t\gamma^{-q}-1)}{t^2}=(t-(\gamma+\gamma^{-1})+t^{-1})(t-(\gamma^q+\gamma^{-q})+t^{-1}).$$ Now taking $q$-th powers, the factors get interchanged, i.e. the above is an element of $\Fq$. Hence $det(s\otimes t-1)$ is a square in $\Fq^*,$ and hence using \ref{character formula}, we see that $\eta(s\otimes t)=1$ for all $t \in SO(E)$.

On the other hand, if $t \in O(E)\backslash SO(E)$, the eigenvalues of $t$ are $\pm1$. Hence $s\otimes t-1$ is invertible with eigenvalues {$\{\gamma-1, \gamma^q-1, \gamma^{-1}-1, \gamma^{-q}-1,-\gamma-1, -\gamma^q-1, -\gamma^{-1}-1, -\gamma^{-q}-1\}.$} Hence $$det(s\otimes t-1)=(1-\gamma^2)(1-\gamma^{2q})(1-\gamma^{-2})(1-\gamma^{-2q})=((\gamma-\gamma^{-1})(\gamma^q-\gamma^{-q}))^2.$$ But now $(\gamma-\gamma^{-1})(\gamma^q-\gamma^{-q}) \notin \Fq$ since it is not fixed by the $q$-th power map for $\gamma \neq \pm1$, hence $\eta(s\otimes t)=-1$ for all $t \in O(E)\backslash SO(E)$. Hence we have proved the lemma.

\epf

Hence substituting the values $\eta(s\otimes t)$ in (\ref{eqn}), we conclude that $\theta_{10}(s) = 1$ and we have \ref{rss}.

\epf

It is now easy to see that $\theta_{10}$ is unipotent.

\bprop
$\theta_{10}$ is unipotent. In fact, we have
\beq
\langle\theta_{10}, R_{T,1}\rangle = 1.
\eeq
\eprop
\begin{proof}
By \ref{fairy tale} and \ref{rss} we see that for $s \in T^F, s\neq \pm1$
\beq
\theta_{10}(s)=1= \sum_{\theta \in \hat{T}^F}{\theta(s)\langle \theta_{10}, R_{T,\theta} \rangle}.
\eeq
Note that $T^F$ has only two non-regular characters, namely the trivial character and the quadratic character $\mu$. It is clear that $\theta_{10}$ is not one of the $\pm R_{T, \theta}$ corresponding to the regular $\theta \in \hat{T}^F$. This is because, for example, the dimension of $R_{T, \theta}$ is $(q^2-1)^2$(using \ref{dimension}), which is not equal to that of $\theta_{10}$.   Hence the regular $\theta$ do not contribute to the sum. By \ref{strong orthogonality} we see that $\theta_{10}$ cannot occur in both $R_{T,1}$ and $R_{T,\mu}$, and by the equation above, cannot occur only in $R_{T,\mu}$. Hence the equation just reads $\langle\theta_{10}, R_{T,1}\rangle = 1$.
\end{proof}

\end{document}